\documentclass[12pt]{article}

\usepackage{pstricks,pst-node}
\usepackage{pathlatex}
\usepackage{psfig}
\usepackage{epsfig}
\usepackage{amsmath,amsthm}
\usepackage{amssymb}
\usepackage{color}
\usepackage{xspace}
\usepackage[colorlinks=true,
linkcolor=green,
filecolor=brown,
citecolor=green]{hyperref}

\def\mbf#1{\mathchoice{\hbox{\boldmath $\displaystyle #1$}}
        {\hbox{\boldmath $\textstyle #1$}}
        {\hbox{\boldmath $\scriptstyle #1$}}
        {\hbox{\boldmath $\scriptscriptstyle #1$}}} % mathboldface

\begin{document}
\newtheorem{theorem}{Theorem}
\newtheorem{defn}[theorem]{Definition}
\newtheorem{lemma}[theorem]{Lemma}
\newtheorem{prop}[theorem]{Proposition}
\newtheorem{cor}[theorem]{Corollary}
%\vspace*{-5mm}
\begin{center}
{\Large
A variant of Touchard's Catalan number identity                         \\ 
}

\vspace{10mm}
David Callan  \\
\noindent {\small Dept. of Statistics, 
University of Wisconsin-Madison,  Madison, WI \ 53706}  \\
{\bf callan@stat.wisc.edu} 

May 12, 2013
\end{center}

\begin{abstract}
It is well known that the Catalan number $C_{n	}$ counts dissections of a regular $(n+2)$-gon into triangles.
Here we count such dissections by number of triangles that contain two sides of the polygon among their three edges, leading to
a combinatorial interpretation of the identity
\[
C_{n} =\sum_{1 \,\le\,k \,\le\, n/2} 2^{n-2k} {n \choose 2k} C_{k}\,\frac{k(n+2)}{n(n-1)},
\]
and illustrating its connection with  Touchard's identity.
\end{abstract}

\section{Introduction}

\vspace*{-5mm}

Consider a regular polygon of $n+2$ sides with one side designated the base. It is a classic result that there are the Catalan number $C_{n}$ ways to insert noncrossing diagonals connecting nonadjacent vertices of the polygon so as to dissect it into triangles (see illustration in Figure \ref{fig}). Each such dissection contains $n-1$ diagonals and $n$ triangles. When $n\ge 2$, each triangle may have 0, 1, or 2 sides in common with the polygon. Let $u_{n,k}$ denote the number of dissections in which precisely $k$ triangles contain 2 sides of the polygon. In any dissection, the number of such 2-polygon-side triangles ranges from a minimum of 2 (provided $n\ge 2$) to a maximum of $\lfloor(n+2)/2\rfloor$.  

Our main result is that $u_{n,k+1}=2^{n-2k} {n \choose 2k} C_{k}\,\frac{k(n+2)}{n(n-1)}$, yielding the apparently new identity
% valid for $n\ge 2$,
\begin{equation}\label{main}
\hspace*{30mm} C_{n} =\sum_{1 \,\le\,k \,\le\, n/2} 2^{n-2k} {n \choose 2k} C_{k}\,\frac{k(n+2)}{n(n-1)},\hspace*{20mm} n\ge 2.
\end{equation}
This identity is reminiscent of Touchard's identity \cite{A091894}, 
\begin{equation}\label{touchard}
\hspace*{38mm} C_{n+1} =\sum_{0 \,\le\,k \,\le\, n/2} 2^{n-2k} {n \choose 2k} C_{k},\hspace*{27mm} n\ge 0,
\end{equation}
%$C_{n+1} =\sum_{0 \,\le\,k \,\le\, n/2} 2^{n-2k} {n \choose 2k} C_{k}$ 
 and indeed we will see a connection between them.
To obtain an expression for $u_{n,k}$, it is convenient to color the base of the polygon blue and the remaining edges black, and let $v_{n,k}$ denote the number of dissections in which $k$ triangles contain two \emph{black} edges of the polygon. In Section 2, we express the \{$u_{n,k}$\} directly in terms of the \{$v_{n,k}$\}. In Section 3, we show bijectively that $v_{n+1,k+1}$ is actually the summand $2^{n-2k} {n \choose 2k} C_{k}$ in (\ref{touchard}), incidentally giving another combinatorial interpretation of Touchard's identity. Section 4 then establishes the main result.

\section[A relation between u(n,k) and v(n,k)]{\protect{A relation between $\mbf{u_{n,k}}$ and $\mbf{v_{n,k}}$}}

\vspace*{-5mm}

Clearly, $u_{1,1}=1,\ u_{2,2}=2,\ u_{3,2}=5$. For $n\ge 4$, let us count the contribution to $u_{n,k}$ according to the positive vertex $r$ of the triangle that contains the base, after labelling the vertices of the polygon $r=-1,0,1,...,n$ counterclockwise from the left endpoint of the base as illustrated in Figure \ref{fig}. We find that the contribution to $u_{n,k}$ for both $r=1$ and $r=n$ is $v_{n-1,k-1}$, and for $2\le r \le n-1$, the contribution is
$
\sum_{k-\frac{n-r+1}{2} \le\, j\, \le \frac{r}{2}}v_{r-1,j}\,v_{n-r,k-j}.
$
Hence,
\[
u_{n,k} = 2 v_{n-1,k-1} + \sum_{r=2}^{n-1}\:\sum_{k-\frac{n-r+1}{2} \le\, j\, \le \frac{r}{2}}v_{r-1,j}\,v_{n-r,k-j},
\]
valid for $n \ge 4 ,\ 2 \le k \le \frac{n + 2}{2}$. Similarly, we find a recurrence for $v_{n,k}$,
\[
v_{n,k} = 2 v_{n-1,k} + \sum_{r=2}^{n-1}\:\sum_{k-\frac{n-r+1}{2} \le\, j\, \le \frac{r}{2}}v_{r-1,j}\,v_{n-r,k-j},
\]
that involves the same double sum. Eliminating the double sum in the two equations leads to the relation
\begin{equation}\label{eq2}
u_{n,k} = v_{n,k} + 2v_{n-1,k-1}-2v_{n-1,k},
\end{equation}
which in fact holds for all $n,k$.

%\textrm{{\normalsize $\frac{n + 2}{2}$}}
\section{A bijection} \label{bij}

\vspace*{-5mm}

It is well known that $2^{n-1-2k} {n-1 \choose 2k} C_{k}$ is the number of Dyck paths that contain $k\ DDU$'s, where $U$ denotes an upstep and $D$ a downstep. (See \cite{dyck2004} for a bijective proof.) We now present a bijection from polygon dissections to Dyck paths which makes it visually obvious that the triangles containing two black sides, taken in clockwise order from the base, except that the last one is ignored, correspond to the $DDU$'s, taken left to right, in the Dyck path. This bijection is simply the composition of the following 3 well known bijections, (1) the  Erdelyi-Etherington bijection from triangle-dissections of a polygon to binary trees \cite[p.\,171]{ec2}, (2) the standard bijection from binary trees to ordered trees (Knuth's ``natural'' correspondence \cite[Section 2.3.2]{acp1}), and (3) the (trivial) ``glove'' bijection from ordered trees to Dyck paths. Here is an illustration.
\setcounter{equation}{0}
%\newpage

\vspace*{-1mm}

\begin{equation}\label{fig}\notag
\textrm{Figure 1}
\end{equation}

\begin{center}

\begin{pspicture}(-1,-4)(10,2.5)%\showgrid
   \psset{unit=1.4cm}  
 
\psdots(-2,0)(-1.62,1.18)(-.62,1.9)(.62,1.9)(1.62,1.18)(2,0)(1.62,-1.18)(.62,-1.9)(-.62,-1.9)(-1.62,-1.18)

\pspolygon(-2,0)(-1.62,1.18)(-.62,1.9)(.62,1.9)(1.62,1.18)(2,0)(1.62,-1.18)(.62,-1.9)(-.62,-1.9)(-1.62,-1.18)(-2,0)

\psline[linecolor=blue,linewidth=1pt](-.62,-1.9)(.62,-1.9)

\psline(-.62,-1.9)(-.62,1.9)
\psline(.62,-1.9)(.62,1.9)
\psline(-.62,1.9)(-2,0)(-.62,-1.9)
\psline(.62,-1.9)(1.62,1.18)(1.62,-1.18)
\psline(-.62,-1.9)(.62,1.9)

\rput(-.8,-2.1){\textrm{{\footnotesize $-1$}}}
\rput(.75,-2.1){\textrm{{\footnotesize $0$}}}
\rput(1.8,-1.18){\textrm{{\footnotesize $1$}}}
\rput(2.2,0){\textrm{{\footnotesize $2$}}}
\rput(-1.8,-1.18){\textrm{{\normalsize $n$}}}
\rput(1.75,1.5){$\ddots$}
\rput(-2.2,0.1){\vdots}

\rput(2.9,-.5){$\longrightarrow$}

\psline[linecolor=red,linewidth=1pt](0,-2.4)(0,-1.4)(-.3,0.5)(-1,0)(-1.5,1)
\psline[linecolor=red,linewidth=1pt](-1,0)(-1.5,-1)
\psline[linecolor=red,linewidth=1pt](0,-1.4)(1,0)(1.3,-.5)(1.8,0)

\rput(0,-2.8){\textrm{{\footnotesize a dissection into $n=8$ triangles}}}
\rput(6,-2.8){\textrm{{\footnotesize left-planted binary tree}}}

\psset{dotscale=1.5} 
\psdots[linecolor=red](0,-2.4)(0,-1.4)(-.3,0.5)(-1,0)(-1.5,1)(-1.5,-1)(1,0)(1.3,-.5)(1.8,0)  

\psset{unit=1.0cm}
\psdots[linecolor=red](5,1)(6,0)(7,1)(7,-1)(8,-2)(9,-3)(9,-1)(10,0)(9,1)
\psline[linecolor=red](5,1)(6,0)(7,1)
\psline[linecolor=red](6,0)(7,-1)(8,-2)(9,-3)
\psline[linecolor=red](8,-2)(9,-1)(10,0)(9,1)

\psset{unit=1.4cm}
\psset{dotscale=1.6} 
\psdots[linecolor=red](-1.5,1)(-1.5,-1) 
\psset{unit=1.0cm}
\psset{dotscale=1.2} 
\psdots[linecolor=red](5,1)(7,1) 

\rput(11.7,-.5){$\xrightarrow[45^{\circ}]{\mathrm{rotate}}$}

\end{pspicture}
\end{center} 
\begin{center} 
\begin{pspicture}(-1,-1)(14,5)%\showgrid
   \psset{unit=1.0cm}  
\psdots(-1,0)(-1,1)(-1,2)(-1,3)(-1,4)(0,3)(0,1)(1,1)(1,2)
\psline(-1,0)(-1,1)(-1,2)(-1,3)(-1,4)
\psline(-1,3)(0,3)
\psline(-1,1)(0,1)(1,1)(1,2)

\rput(3,2){$\xrightarrow[\textrm{east$\:\rightarrow\:$west edges}]{\textrm{swing down}}$}

\psdots(5,0)(5,1)(5,2)(5,3)(5,4)(6,3)(6,1)(7,1)(7,2)
\psline(6,3)(5,2)(5,3)(5,4)
\psline(6,1)(5,0)(5,1)(5,2)
\psline(5,0)(7,1)(7,2)

\rput(8.2,2){$\xrightarrow{\textrm{``prettify''}}$}

\psdots(9.5,4)(9.5,3)(10,2)(10,1)(10.5,3)(11,0)(11,1)(12,1)(12,2)
\psline(9.5,4)(9.5,3)(10,2)(10.5,3)
\psline(10,2)(10,1)(11,0)(11,1)
\psline(11,0)(12,1)(12,2)
\rput(11,-0.5){\textrm{{\footnotesize ordered tree}}}

\rput(13,2){$\longrightarrow$}
\psset{dotscale=1.8} 
\psdots(-1,4)(0,3)(5,4)(6,3)(9.5,4)(10.5,3)

\end{pspicture}

\end{center}
\vspace*{-3mm}
\Einheit=0.6cm
\[
\Pfad(-8,0),3333443444343344\endPfad
\SPfad(-8,0),1111111111111111\endSPfad
\DuennPunkt(-8,0)
\DuennPunkt(-7,1)
\DuennPunkt(-6,2)
\DuennPunkt(-5,3)
\NormalPunkt(-4,4)
\DuennPunkt(-3,3)
\DuennPunkt(-2,2)
\NormalPunkt(-1,3)
\DuennPunkt(0,2)
\DuennPunkt(1,1)
\DuennPunkt(2,0)
\DuennPunkt(3,1)
\DuennPunkt(4,0)
\DuennPunkt(5,1)
\DuennPunkt(6,2)
\DuennPunkt(7,1)
\DuennPunkt(8,0)
%\Label\u{\textrm{ {\scriptsize Dyck path}}}(0,-.5)
\] 
\begin{center} 
\begin{pspicture}(-1,-1)(0,0)%\showgrid
   \psset{unit=1.0cm}  
\rput(-.5,0.3){\textrm{{\footnotesize Dyck path}}} 
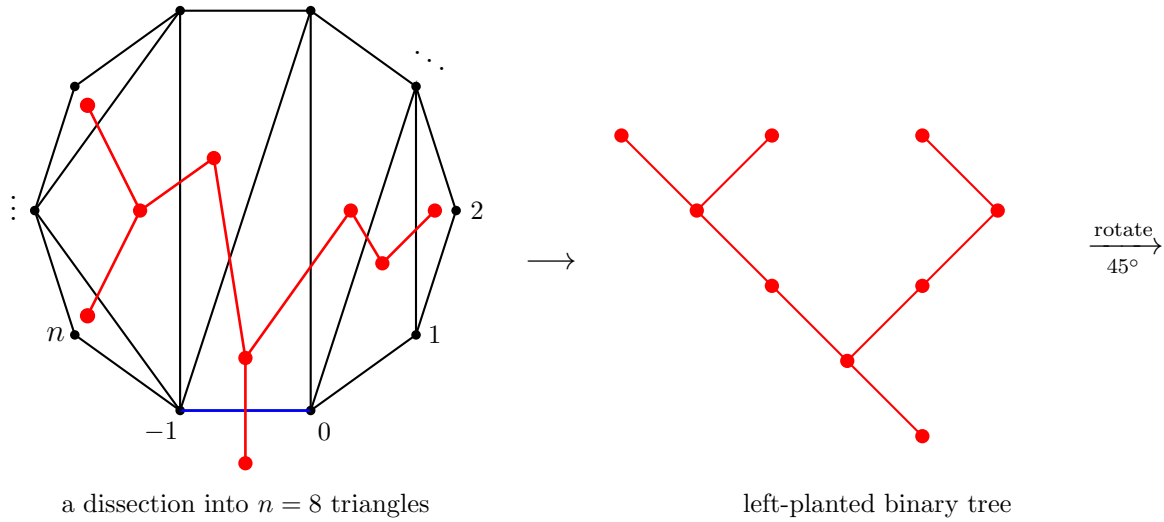
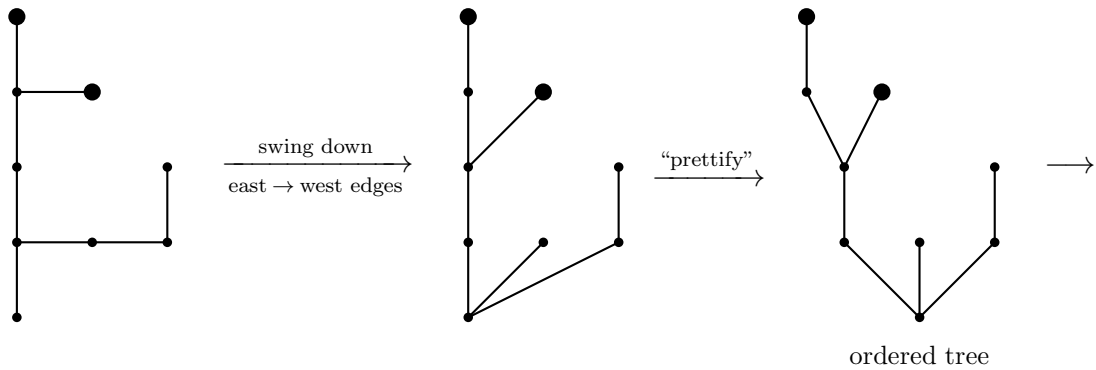
\rput(0,-.5){\textrm{ {\normalsize Bijection from triangle dissections to Dyck paths}}}
\end{pspicture}
\end{center}

%\vspace*{-15mm}

%\begin{table}\label{fig}\notag
%\textrm{Figure 1}
%\end{table}

\vspace*{-5mm}

\noindent The last step is the glove bijection: walk clockwise around the tree starting from the root and record an upstep (resp. downstep) each time an edge is traversed upward (resp. downward). Or, more picturesquely, burrow up the edges from the root to form a multi-fingered glove and fan out the fingers. Thus each edge in the tree corresponds to a matching upstep and downstep in the path. 

The illustrated dissection has 3 triangles containing two black sides of the polygon; all but the last are highlighted using enlarged dots, and they show up in the Dyck path as vertices initiating  a descent of 2 or more downsteps followed by an upstep, that is, they correspond to $DDU$s in the Dyck path, as claimed.

\vspace*{-2mm}

\section{Conclusion}

\vspace*{-5mm}

The preceding section shows that $v_{n,k} = 2^{n+1-2k}\binom{n-1}{2k-2}C_{k-1}$. Substituting into (\ref{eq2}), we find 
\[
u_{n,k}=2^{n+1-2k}\left(\binom{n-2}{2k-3}C_{k-1}+\binom{n-2}{2k-4}4C_{k-2}\right),
\]
which simplifies to 
\[
u_{n,k+1}=2^{n-2k} \binom{n}{2k}C_{k} \frac{(n+2)k}{n(n-1)}.
\]
Now sum over $k$ to obtain (\ref{main}). The first few values of $u_{n,k}$ are given in the following table.

\[
\begin{array}{c|ccccc}
	n^{\textstyle{\,\backslash \,k}} &  2 & 3 & 4 & 5  \\
\hline 	 
 	2&    2 & & &   \\
	3&    5 &  & &   \\
	4&    12 & 2  &   &    \\ 
	5&    28 & 14 &   &    \\
	6&    64 & 64 & 4 &    \\
	7&    144 & 240 & 45 &   \\
    8&    320 & 800 & 300 & 10    \\
    9&    704 & 2464 & 1540 & 154   \\ 

\end{array}
\]

\vspace*{2mm}

\centerline{{\normalsize Table of values of $u_{n,k}$ }}

\vspace*{4mm}

\textbf{Added in proof} \ Tewodros Amdeberhan informs me that he has recently discovered an identity equivalent to (\ref{main}), namely
\[
\frac{2n}{n+3} C_{n+1} =\sum_{0 \,\le\,k \,\le\, (n-1)/2} 2^{n-2k} {n \choose 2k+1} C_{k}\,\frac{2k+1}{k+2},
\]
and has observed that subtracting the latter from Touchard's identity (\ref{touchard})(multiplied by 2) gives an alternating sum expression for the super ballot number $6/(n+3)C_{n+1}$, sequence \htmladdnormallink{A007054}{http://oeis.org/A007054} in OEIS.

\textbf{Acknowledgement of priority}\ \ Alon Regev pointed out to me that the main results of this paper have previously been obtained by Hurtado and Noy \cite{ears}.

\end{document}